\documentclass[12pt]{aptpub}

\renewcommand{\crnotice}{\hfill{{\small\em

\/}\ {\small}}}

\usepackage{amsmath,amstext,url}
\usepackage{color}
\usepackage{amsfonts}
\usepackage{amsmath}
\usepackage{txfonts}
\usepackage{mathrsfs}
\usepackage{cite}

\renewcommand{\authornames}{\underline{\footnotesize }
J. Li, W. Zhang}
\renewcommand{\shorttitle}
{ \underline{{The multiple birth properties of multi-type Markov branching processes}}}

\numberwithin{equation}{section} 

\makeatletter \@addtoreset{equation}{section}

\makeatletter \@addtoreset{lemma}{section}

\makeatletter \@addtoreset{theorem}{section}

\makeatletter \@addtoreset{proposition}{section}

\makeatletter \@addtoreset{corollary}{section}

\makeatletter \@addtoreset{remark}{section}

\makeatletter \@addtoreset{definition}{section}

\makeatletter \@addtoreset{example}{section}

\begin{document}

\thispagestyle{firstpg}

\vspace*{1.5pc} \noindent \normalsize\textbf{\Large {The multiple birth properties of multi-type Markov branching processes}} \hfill

\vspace{12pt} \hspace*{0.75pc}{\small\textrm{\uppercase{Junping Li}}}
\hspace{-2pt}$^{*}$, {\small\textit{Guangdong University of Science $\&$ Technology; Central
South University}

\hspace*{0.75pc}{\small\textrm{\uppercase{Wanting Zhang}}}\hspace{-2pt}$^{**}$, {\small\textit{Central South University }}}


\par
\footnote{\hspace*{-0.75pc}$^{*}\,$Postal address:
 Guangdong University of Science $\&$ Technology, Dongguan, 523083, China; Central
South University, Changsha, 410083, China. E-mail:
jpli@mail.csu.edu.cn}

\par
\footnote{\hspace*{-0.75pc}$^{**}\,$Postal
address: School of Mathematics and Statistics, Central
South University, Changsha, 410083, China. E-mail address:
202111062@csu.edu.cn}
\par
\renewenvironment{abstract}{%
\vspace{8pt} \vspace{0.1pc} \hspace*{0.25pc}
\begin{minipage}{14cm}
\footnotesize
{\bf Abstract}\\[1ex]
\hspace*{0.5pc}} {\end{minipage}}
\begin{abstract}
   The main purpose of this paper is to consider the multiple birth properties for multi-type Markov branching processes. We first construct a new multi-dimensional Markov process based on the multi-type Markov branching process, which can reveal the multiple birth characteristics. Then the joint probability distribution of multiple birth of multi-type Markov branching process until any time $t$ is obtained by using the new process. Furthermore, the probability distribution of multiple birth until the extinction of the process is also given.
\end{abstract}
\vspace*{12pt} \hspace*{2.25pc}
\parbox[b]{26.75pc}{{
}}
\par
\normalsize
\par
\renewcommand{\amsprimary}[1]{
     \vspace*{8pt}
     \hspace*{2.25pc}
     \parbox[b]{12.75pc}
{\scriptsize
AMS 2000 Subject Classification: Primary 60J27 Secondary 60J35
     {\uppercase{#1}}}\par\normalsize}
\renewcommand{\ams}[2]{
     \vspace*{8pt}
     \hspace*{2.25pc}
     \parbox[b]{18.75pc}{\scriptsize
     AMS 2000 SUBJECT CLASSIFICATION: PRIMARY
     {\uppercase{#1}}\\ \phantom{
     AMS 2000
     SUBJECT CLASSIFICATION:
     }
    SECONDARY
 {\uppercase{#2}}}\par\normalsize}

\ams{60J27}{60J35}

\par
\vspace{5mm}
 \setcounter{section}{1}
 \setcounter{equation}{0}
 \setcounter{theorem}{0}
 \setcounter{lemma}{0}
 \setcounter{corollary}{0}
\noindent {\large \bf 1. Introduction}
\vspace{3mm}
\par
Markov branching processes play an important role in the research and applications of stochastic processes. Standard references are Anderson~\cite{And91}, Harris~\cite{Har63}, Athreya \& Ney~\cite{Ath-N72}, Asmussen \& Hering~\cite{Asm-H83}, Athreya \& Jagers~\cite{Asm-J97} and others.
\par
The basic property governing the evolution of a Markov branching process is the branching property, i.e., different individuals act independently when giving offsprings. The classical Markov branching processes are well studied, some related references are Harris~\cite{Har63}, Athreya \& Ney~\cite{Ath-N72}, Asmussen \& Hering~\cite{Asm-H83}, and Athreya \& Jagers~\cite{Asm-J97}. Based on the branching structure, there are many references concentrating on generalization of ordinary Markov branching processes. For example, Vatutin~\cite{Vat74}, Li, Chen \& Pakes~\cite{Li-C-P12} considered the branching processes with state-independent immigration. Chen, Li \& Ramesh~\cite{Chen-L-R05} and Chen, Pollet, Zhang \& Li~\cite{Chen-P-L-Z07} considered weighted Markov branching processes, Li \& Chen~\cite{Li-C18} considered generalized Markov interacting branching processes, Li \& Wang~\cite{Li09,Li-W12} and Meng \& Li~\cite{Meng-L18} considered $n$-type branching processes with or without immigration. Recently, Li \& Li~\cite{Li-Li21,Li-Li22} considered down/up crossing properties of weighted Markov collision processes and one-dimensional Markov branching processes.
\par
In this paper, we mainly discuss the multiple birth properties of multi-type Markov branching processes. Different from the one-type case, the number of individuals of other types may change when an individual splits.
\par
For convenience of our discussion, we make the following notations throughout of this paper. Let $\mathbf{Z}_+$ be the set of nonnegative integers.
\par
(C-1)\ $\mathbf{Z}_+ ^d :=\{\emph{\textbf{i}}=({i_1},\cdots,{i_d}):{i_1},\cdots,{i_d} \in {\mathbf{Z}_+}\}$, and for any $\emph{\textbf{i}}=(i_1,\cdots,i_d) \in \mathbf{Z}_+^d$, denote $\mid\emph{\textbf{i}}\mid=\sum\limits_{k=1}^d i_k$.
\par
(C-2)\ ${[0,1]^d} = \{\emph{\textbf{x}}=({x_1},\cdots,{x_d}):0 \leq {x_1},\cdots,{x_d} \leq 1\} $.
\par
(C-3)\ ${\chi_{_{\mathbf{Z}_+^d}}}(\cdot )$ is the indicator of $\mathbf{Z}_+^d$
\par
(C-4)\ $\emph{\textbf{0}}=(0,\cdots,0)$, $\emph{\textbf{1}}= (1,\cdots,1)$, ${\emph{\textbf{e}}_k} = (0,\cdots,{1_k},\cdots,0)$ are vectors in ${[0,1]^d}$.
\par
(C-5)\ For any $\emph{\textbf{x}},\emph{\textbf{y}}\in [0,1]^d$, $\emph{\textbf{x}}\leq \emph{\textbf{y}}$ means $x_k \leq y_k$ for all $k= 1,\cdots,d$. $\emph{\textbf{x}}<\emph{\textbf{y}}$ means $x_k \leq y_k$ for all $k= 1,\cdots,d$, and $x_k<y_k$ for at least one $k$.
\par
(C-6)\ For any $\emph{\textbf{x}}\in [0,1]^d$, denote $\| \emph{\textbf{x}}\|_1 = \sum\limits_{k=1}^d |x_k|$.
\par
A $d$-type Markov branching process can be intuitively described as follows:
\par
(1)\ Consider a system involving $d$ types of individuals. The life length of a type-$k$ individual is exponentially distributed with mean ${\mathbf{\theta}_k}\ (k= 1,\cdots,d)$.
\par
(2)\ Individuals in the system split independently. When a type-$k$ individual dies after a random time, it is replaced by ${j_1}$ individuals of type-$1$, $\cdots$, and $j_d$ individuals of type-$d$, with probability $p^{(a)_{\emph{\textbf{j}}}}$, here $\emph{\textbf{j}}=(j_1,\cdots,j_d)$. Without loss of generality, we can assume $p^{(k)}_{\emph{\textbf{e}}_k}=0\ (k=1,\cdots,d)$ since such split does not change the state of the system.
\par
(3)\ When this system is empty, it stops. i.e., $\emph{\textbf{0}}$ is an absorbing state.
\par
We now define the infinitesimal generator of $d$-type Markov branching processes, i.e., the $Q$-matrix.
\par
\begin{definition}\label{def1.1}\
 A $Q$-matrix $Q = (q_{\emph{\textbf{ij}}}:\emph{\textbf{i}},\emph{\textbf{j}} \in \mathbf{Z}_+^d)$ is called a $d$-type Markov branching $Q$-matrix (henceforth referred to as a $d$TMB $Q$-matrix), if
\begin{eqnarray}\label{eq1.1}
q_{\emph{\textbf{ij}}}=\begin{cases}
\sum\limits_{k=1}^d i_kb_{\emph{\textbf{j}}-\emph{\textbf{i}}- \emph{\textbf{e}}_k}^{(k)},\ & if \mid\emph{\textbf{i}}\mid> 0,\\
0,\ &\ otherwise.
\end{cases}
\end{eqnarray}
where $b^{(k)}_{\emph{\textbf{j}}}=0$ for $\emph{\textbf{j}}\notin \mathbf{Z}_+^d$ and
\begin{eqnarray}\label{eq1.2}
b_{\emph{\textbf{j}}}^{(k)}=\theta_kp_{\emph{\textbf{j}}}^{(k)} \geq 0\ (\ \emph{\textbf{j}}\neq \emph{\textbf{e}}_k), \quad b^{(k)}_{\emph{\textbf{e}}_k}=-\sum\limits_{\emph{\textbf{j}}\neq \emph{\textbf{e}}_k} b_{\emph{\textbf{j}}}^{(k)}\ (k=1,\cdots,d).
\end{eqnarray}
\end{definition}
\par
\begin{definition}\label{def1.2}\
A $d$-type Markov branching process (henceforth referred to as $d$TMBP) is a continuous time Markov chain with state space $\mathbf{Z}_ + ^d$ whose transition probability function $P(t)=(p_{\emph{\textbf{ij}}}(t):\emph{\textbf{i}},\emph{\textbf{j}}\in \mathbf{Z}_+^d)$ satisfies the Kolmogorov forward equation
\begin{eqnarray*}
    P'(t)=P(t)Q.
\end{eqnarray*}
where $Q$ is given in (\ref{eq1.1})-(\ref{eq1.2}),
\end{definition}
\par
\vspace{5mm}
\setcounter{section}{2}
\setcounter{equation}{0}
\setcounter{theorem}{0}
\setcounter{lemma}{0}
\setcounter{definition}{0}
\setcounter{corollary}{0}
\noindent {\large \bf 2. Preliminaries}
\vspace{3mm}
\par
In this section, we make some preliminaries related to the problem considered in this paper. For $k=1,\cdots, d$, let $R_a\subset \mathbf{Z}_+^d$ be finite subsets with $b^{(k)}_{\emph{\textbf{j}}}> 0$ for any $\emph{\textbf{j}}\in R_k$. Also let $r_k$ denote the number of elements in $R_k$ and $r=r_1+\cdots+r_d$. This paper is devoted to considering the probability distribution property of the number of type-$k$ individuals giving $R_k$-birth until time $t$.
\par
For convenience of our discussion, we only discuss the case of $2$-type Markov branching process. The general case of the $d$-type $(d \geq 3)$ can be studied analogously.
\par
Define
\begin{eqnarray}\label{eq2.1}
B_k(\emph{\textbf{x}})=\sum\limits_{\emph{\textbf{j}}\in {\mathbf{Z}}_+^2}b^{(k)}_{\emph{\textbf{j}}}\emph{\textbf{x}}
^{\emph{\textbf{j}}},\quad \emph{\textbf{x}}\in [0,1]^2,\ \ k=1,2,
\end{eqnarray}
and
\begin{eqnarray*}
B_{ij}(\emph{\textbf{x}})=\frac{\partial B_i(\emph{\textbf{x}})}{\partial x_j},\quad \emph{\textbf{x}}\in [0,1]^2,\ \ i,j=1,2.
\end{eqnarray*}
\par
In order to avoid some trivial cases, we assume the following conditions hold.
\par
(A-1)\ $(B_1(\emph{\textbf{x}}),B_2(\emph{\textbf{x}}))$ is nonsingular, i.e., there is no $2\times 2$-matrix $M$ such that $(B_1(\emph{\textbf{x}}),B_2(\emph{\textbf{x}}))=\emph{\textbf{x}}M$;
\par
(A-2)\ $B_{ij}(1,1)<\infty,\ i,j=1,2$;
\par
(A-3)\ The matrix $(B_{ij}(1,1):i,j=1,2)$ is positively regular, i.e., there exists an integer $m$ such that $(B_{ij}(1,1):i,j=1,2)^m>0$ in sense of all the elements are positive.
\par
For any $\emph{\textbf{x}}\in [0,1]^2$, the maximal eigenvalue of $(B_{ij}(1,1):i,j=1,2)$ is denoted by $\rho(\emph{\textbf{x}})$. The following lemma is due to Li \& Wang~\cite{Li-W12}, we only state it without proof.
\par
\begin{lemma}\label{le2.1}
The system of equations
\begin{eqnarray}\label{eq2.a}
 \begin{cases} B_1(\textbf{x})=0, \\
               B_2(\textbf{x})=0.
 \end{cases}
\end{eqnarray}
has at most two solutions in $[0,1]^2$. Let $\textbf{q}=(q_1,q_2)$ denote the smallest
nonnegative solution to $(\ref{eq2.a})$. Then,
\par
{\rm (i)}\ $q_i$ is the extinction probability when the Feller
minimal process starts at state $\textbf{e}_i\ (i=1,2)$.
Moreover, if $\rho(\textbf{1})\leq 0$, then $\textbf{q}=\textbf{1}$; while
if $\rho(\textbf{1})>0$, then $\textbf{q}<\textbf{1}$, i.e.,
$q_1, q_2<1$.
\par
{\rm (ii)}\  $\rho(\textbf{q})\leq 0$.
\end{lemma}
\par
The following result is well-known which reveals the basic property of $2$-type Markov branching processes.
\par
\begin{lemma}
Let $P(t)=(p_{\textbf{i}\textbf{j}}(t):\textbf{i},\textbf{j}\in \mathbf{Z}_+^2)$ be the transition function with $Q$-matrix $Q$ given in {\rm(\ref{eq1.1})-(\ref{eq1.2})}. Then,
\begin{eqnarray*}
\frac{\partial F_{\textbf{i}}(t,\textbf{x})}{\partial t}
=B_1(\textbf{x})\frac{\partial F_{\textbf{i}}(t,\textbf{x})}{\partial x_1}+B_2(\textbf{x})\frac{\partial F_{\textbf{i}}(t,\textbf{x})}{\partial x_2},
\end{eqnarray*}
where $F_{\textbf{i}}(t,\textbf{x})=\sum\limits_{\textbf{j}\in \mathbf{Z}_+^2}p_{\textbf{i}\textbf{j}}(t)\textbf{x}
^{\textbf{j}}$ with $\textbf{x}^{\textbf{j}}=x_1^{j_1}x_2^{j_2}$.
\end{lemma}
\par
Li \& Meng~\cite{Li-M17} derived the regularity criteria for $2$-type Markov branching processes. Assumption (A-1) guarantees the regularity of the process.
\par
Let $\emph{\textbf{Y}}(t)=(Y_{\emph{\textbf{k}}}(t):\emph{\textbf{k}}\in R_1)$ be the number of type-$1$ individuals giving $R_1$-birth until time $t$ and $\emph{\textbf{Z}}(t)=(Z_{\emph{\textbf{k}}}(t):\emph{\textbf{k}}\in R_2)$ be the number of type-$2$ individuals giving $R_2$-birth until time $t$. We will discuss the probability distribution property of $(\emph{\textbf{Y}}(t),\emph{\textbf{Z}}(t))$. For this end, we define
\begin{eqnarray}\label{eq2.2}
B_1(\emph{\textbf{x}},\emph{\textbf{y}}) = \sum\limits_{\emph{\textbf{j}}\in R_1} b_{\emph{\textbf{j}}}^{(1)}\emph{\textbf{x}}^{\emph{\textbf{j}}}
y_{\emph{\textbf{j}}} ,\quad \bar{B}_1(\emph{\textbf{x}})= \sum\limits_{\emph{\textbf{j}}\in R_1^c} b_{\emph{\textbf{j}}}^{(1)}\emph{\textbf{x}}^{\emph{\textbf{j}}}.
\end{eqnarray}
\begin{eqnarray}\label{eq2.3}
B_2(\emph{\textbf{x}},\emph{\textbf{z}}) = \sum\limits_{\emph{\textbf{j}}\in R_2} b_{\emph{\textbf{j}}}^{(2)}\emph{\textbf{x}}^{\emph{\textbf{j}}}
z_{\emph{\textbf{j}}} ,\quad \bar{B}_2(\emph{\textbf{x}})
=\sum\limits_{\emph{\textbf{j}}\in R_2^c} b_{\emph{\textbf{j}}}^{(2)}\emph{\textbf{x}}^{\emph{\textbf{j}}}.
\end{eqnarray}
where $\emph{\textbf{x}}=(x_1,x_2)\in \mathbf{Z}_+^2$; $\emph{\textbf{y}}=(y_{\emph{\textbf{j}}}:\emph{\textbf{j}}\in R_1), \emph{\textbf{z}}=(z_{\emph{\textbf{j}}}:\emph{\textbf{j}}\in R_2)$. It is obvious that $\bar{B}_1(\emph{\textbf{x}})$, $\bar{B}_2(\emph{\textbf{x}})$ are well defined at least on ${[0,1]^2}$. $B_1(\emph{\textbf{x}},\emph{\textbf{y}})$, $B_2(\emph{\textbf{x}},\emph{\textbf{z}})$ are well defined at least on $[0,1]^{2+r_1}$ and $[0,1]^{2+r_2}$ respectively.
\par
Since the $2$-type branching process itself can not to reveal the detailed multi-birth directly, we define a new $Q$-matrix $\tilde{Q} =(q_{_{(\emph{\textbf{i}},\emph{\textbf{k}},\tilde{\emph{\textbf{k}}}),
(\emph{\textbf{j}},\emph{\textbf{l}},\tilde{\emph{\textbf{l}}})}}:
(\emph{\textbf{i}},\emph{\textbf{k}},\tilde{\emph{\textbf{k}}}),
(\emph{\textbf{j}},\emph{\textbf{l}},\tilde{\emph{\textbf{l}}})\in \mathbf{Z}_+^{2+r_1+r_2})$ as follows:
\begin{eqnarray}
 q_{_{(\emph{\textbf{i}},\emph{\textbf{k}},\tilde{\emph{\textbf{k}}}),
(\emph{\textbf{j}},\emph{\textbf{l}},\tilde{\emph{\textbf{l}}})}}=
\begin{cases}
\sum\limits_{a=1}^2 i_ab^{(a)}_{\emph{\textbf{j}}-\emph{\textbf{i}}+ \emph{\textbf{e}}_a},& if\
|~\emph{\textbf{i}}~|>0,\
\emph{\textbf{l}}=\emph{\textbf{k}}+I_{_{R_1}}(\emph{\textbf{j}}
\!-\!\emph{\textbf{i}}\!+\!\emph{\textbf{e}}_1)\varepsilon_{\emph{\textbf{j}}
\!-\emph{\textbf{i}}\!+\emph{\textbf{e}}_1},\ \tilde{\emph{\textbf{l}}}=\tilde{\emph{\textbf{k}}}+I_{_{R_2}}
(\emph{\textbf{j}}
\!-\!\emph{\textbf{i}}\!+\!\emph{\textbf{e}}_2)\tilde{\varepsilon}
_{\emph{\textbf{j}}
\!-\emph{\textbf{i}}\!+\emph{\textbf{e}}_2},\\
0,& otherwise,
\end{cases}
\end{eqnarray}
where $\varepsilon_{\emph{\textbf{k}}}\ (\emph{\textbf{k}}\in R_1)$ denotes the vector in $\mathbf{Z}_+^{r_1}$ with the $\emph{\textbf{k}}$'th element being $1$ and the others being $0$. $\tilde{\varepsilon}_{\tilde{\emph{\textbf{k}}}}\ (\tilde{\emph{\textbf{k}}}\in R_2)$ denotes the vector in $\mathbf{Z}_+^{r_2}$  with the $\tilde{\emph{\textbf{k}}}$'th element being $1$ and the others being $0$. $I_{_{R_1}}$ and $I_{_{R_2}}$ are the indicators of $R_1$ and $R_2$ respectively.
\par
It is obvious that $\tilde{Q}$ determines a $(2+r_1+r_2)$-dimensional continuous-time Markov chain $(\emph{\textbf{X}}(t),\emph{\textbf{Y}}(t),\emph{\textbf{Z}}(t))$, where $\emph{\textbf{X}}(t)$ is the $2$-type Markov branching process, $\emph{\textbf{Y}}(t) =(Y_{\emph{\textbf{k}}}(t):\emph{\textbf{k}} \in R_1)$ (or $\emph{\textbf{Z}}(t) =(Z_{\emph{\textbf{k}}}(t):\emph{\textbf{k}} \in R_2)$) counts the number of type-$1$ (or type-$2$) individuals giving $R_1$-birth (or $R_2$-birth) until time $t$. We assume that $Y_{\emph{\textbf{k}}}(0)=0$ and $Z_{\emph{\textbf{k}}}(0)=0$ for all $\emph{\textbf{k}}\in R_1$ and $\emph{\textbf{k}}\in R_2$. In particular,
\par
(1)\ if $R_1=\{\emph{\textbf{0}}\}$ (or $R_2=\{\emph{\textbf{0}}\}$), then $Y_{\emph{\textbf{0}}}(t)$ (or $Z_{\emph{\textbf{0}}}(t)$) counts the pure death number of type-$1$ (or type-$2$) individuals until time $t$.
\par
(2)\ If $R_1=\{(n_1,n_2)\}$, then $Y_{(n_1,n_2)}(t)$ counts the $(n_1,n_2)$-birth number of type-$1$ individuals until time $t$.
\par
(3)\ If $R_2=\{(n_1,n_2)\}$, then $Z_{(n_1,n_2)}(t)$ counts the $(n_1,n_2)$-birth number of type-$2$ individuals until time $t$.
\par
Let $\tilde{P}(t):=(\tilde{p}_{_{(\emph{\textbf{i}},\emph{\textbf{k}},
\tilde{\emph{\textbf{k}}}),
(\emph{\textbf{j}},\emph{\textbf{l}},\tilde{\emph{\textbf{l}}})}}(t):
(\emph{\textbf{i}},\emph{\textbf{k}},\tilde{\emph{\textbf{k}}}),
(\emph{\textbf{j}},\emph{\textbf{l}},\tilde{\emph{\textbf{l}}})\in \mathbf{Z}_+^{2+r_1+r_2})$ be the transition probability of $(\emph{\textbf{X}}(t),\emph{\textbf{Y}}(t),\emph{\textbf{Z}}(t))$.
Define
\begin{equation*}
 F_{_{\emph{\textbf{i}},\emph{\textbf{k}},
\tilde{\emph{\textbf{k}}}}}(t,\emph{\textbf{x}},\emph{\textbf{y}},
\emph{\textbf{z}})
=\sum\limits_{(\emph{\textbf{j}},\emph{\textbf{l}},
\tilde{\emph{\textbf{l}}})\in \mathbf{Z}_+^{2+r_1+r_2}}\tilde{p}_{_{(\emph{\textbf{i}},\emph{\textbf{k}},
\tilde{\emph{\textbf{k}}}),
(\emph{\textbf{j}},\emph{\textbf{l}},\tilde{\emph{\textbf{l}}})}}(t)
\emph{\textbf{x}}^{\emph{\textbf{j}}}\emph{\textbf{y}}^{\emph{\textbf{l}}}
\emph{\textbf{z}}^{\tilde{\emph{\textbf{l}}}},\quad (\emph{\textbf{x}},\emph{\textbf{y}},\emph{\textbf{z}})\in [0,1]^{2+r_1+r_2},
\end{equation*}
where $\emph{\textbf{x}}^{\emph{\textbf{j}}}=x_1^{j_1}x_2^{j_2}$, $ \emph{\textbf{y}}^{\emph{\textbf{l}}}=\prod\limits_{\emph{\textbf{m}}\in R_1}y_{\emph{\textbf{m}}}^{l_{\emph{\textbf{m}}}}$ and $ \emph{\textbf{z}}^{\tilde{\emph{\textbf{l}}}}=\prod\limits_{\emph{\textbf{m}}\in R_2}z_{\emph{\textbf{m}}}^{\tilde{l}_{\emph{\textbf{m}}}}$.
\par
\begin{lemma}\label{le2.3}
Let $\tilde{P}(t)=(\tilde{p}_{_{(\textbf{i},\textbf{k},
\tilde{\textbf{k}}),
(\textbf{j},\textbf{l},\tilde{\textbf{l}})}}(t):
(\textbf{i},\textbf{k},\tilde{\textbf{k}}),
(\textbf{j},\textbf{l},\tilde{\textbf{l}})\in \mathbf{Z}_+^{2+r_1+r_2})$ be the transition probability of $(\textbf{X}(t),\textbf{Y}(t),\textbf{Z}(t))$. Then,
\par
{\rm(1)}\ for any $(\textbf{x},\textbf{y},\textbf{z}) \in [0,1]^{2+r_1+r_2}$,
\begin{eqnarray}\label{eq2.5}
&&\frac{\partial F_{\textbf{i},\textbf{0},\tilde{\textbf{0}}}(t,\textbf{x},\textbf{y},
\textbf{z})}
{\partial t}\nonumber\\
&=&[B_1(\textbf{x},\textbf{y})+\bar{B}
_1(\textbf{x})]\frac{\partial F_{\textbf{i},\textbf{0},\tilde{\textbf{0}}}(t,\textbf{x},\textbf{y},
\textbf{z})}{{\partial {x_1}}}+[B_2(\textbf{x},\textbf{z})+\bar{B}
_2(\textbf{x})]\frac{\partial F_{\textbf{i},\textbf{0},\tilde{\textbf{0}}}(t,\textbf{x},\textbf{y},
\textbf{z})}{{\partial x_2}}
\end{eqnarray}
where $B_1(\textbf{x},\textbf{y})$, $B_2(\textbf{x},\textbf{z})$, $\bar{B}
_1(\textbf{x})$ and $\bar{B}
_2(\textbf{x})$ are defined in {\rm(\ref{eq2.1})-(\ref{eq2.3})}.
\par
{\rm(2)}\ For any $(\textbf{x},\textbf{y},\textbf{z}) \in [0,1]^{2+r_1+r_2}$ and $(\textbf{i},\textbf{k},\tilde{\textbf{k}}) \in \mathbf{Z}_+^{2+r_1+r_2}$,
\begin{eqnarray}\label{eq2.6}
F_{_{\textbf{i},\textbf{k},
\tilde{\textbf{k}}}}(t,\textbf{x},\textbf{y},
\textbf{z})
=\textbf{y}^{\textbf{k}}\textbf{z}^{\tilde{\textbf{k}}} [\textbf{F}(t,\textbf{x},\textbf{y},
\textbf{z})]^{\textbf{i}}
\end{eqnarray}
where $\textbf{F}(t,\textbf{x},\textbf{y},
\textbf{z})=(F_1(t,\textbf{x},\textbf{y},
\textbf{z}),F_2(t,\textbf{x},\textbf{y},
\textbf{z}))$ with $F_k(t,\textbf{x},\textbf{y},
\textbf{z})=F_{\textbf{e}_k,\textbf{0},\textbf{0}}(t,\textbf{x},\textbf{y},
\textbf{z})\ (k=1,2)$.
\end{lemma}
\par
\begin{proof}
(1)\ By the Kolmogorov forward equation, for any $(\emph{\textbf{i}},\emph{\textbf{k}},
\tilde{\emph{\textbf{k}}}),(\emph{\textbf{j}},\emph{\textbf{l}},
\tilde{\emph{\textbf{l}}})\in \mathbf{Z}_+^{2+r_1+r_2}$,
\begin{equation*}
\tilde{p}'_{(\emph{\textbf{i}},\emph{\textbf{k}},
\tilde{\emph{\textbf{k}}}),
(\emph{\textbf{j}},\emph{\textbf{l}},\tilde{\emph{\textbf{l}}})}(t)
=\sum\limits_{(\emph{\textbf{a}},\emph{\textbf{m}},\tilde{\emph{\textbf{m}}})
\in \mathbf{Z}_+^{2+r_1+r_2}}\tilde{p}_{(\emph{\textbf{i}},\emph{\textbf{k}},
\tilde{\emph{\textbf{k}}}),
(\emph{\textbf{a}},\emph{\textbf{m}},\tilde{\emph{\textbf{m}}})}(t)
q_{(\emph{\textbf{a}},\emph{\textbf{m}},\tilde{\emph{\textbf{m}}}),
(\emph{\textbf{j}},\emph{\textbf{l}},\tilde{\emph{\textbf{l}}})}.
\end{equation*}
Multiplying $\emph{\textbf{x}}^{\emph{\textbf{j}}}\emph{\textbf{y}}^{\emph{\textbf{l}}}
\emph{\textbf{z}}^{\tilde{\emph{\textbf{l}}}}$ on both sides of the above equation and summing over $(\emph{\textbf{j}},\emph{\textbf{l}},\tilde{\emph{\textbf{l}}})\in \mathbf{Z}_+^{2+r_1+r_2}$ yield (\ref{eq2.5}).
\par
(2)\ Let $\emph{\textbf{X}}_{a,k}(t)$ denote the offsprings at time $t$ of the $k$'th individual of type-$a$ at initial, $\emph{\textbf{Y}}_{a,k}(t)$ denote the number of $R_1$-birth individuals of $\emph{\textbf{X}}_{a,k}(t)\ (a=1,2)$ and $\emph{\textbf{Z}}_{a,k}(t)$ denote the number of $R_2$-birth individuals of $\emph{\textbf{X}}_{a,k}(t)\ (a=1,2)$. Then, $\{(\emph{\textbf{X}}_{a,k}(t),\emph{\textbf{Y}}_{a,k}(t),
\emph{\textbf{Z}}_{a,k}(t)): k=1,\cdots, i_a; a=1,2\}$ are independent. Moreover, for $a=1,2$, $(\emph{\textbf{X}}_{a,k}(t),\emph{\textbf{Y}}_{a,k}(t),
\emph{\textbf{Z}}_{a,k}(t))$ has the common distribution of $(\emph{\textbf{X}}(t),\emph{\textbf{Y}}(t),
\emph{\textbf{Z}}(t))$ starting at $(\emph{\textbf{e}}_a,\emph{\textbf{0}},\emph{\textbf{0}})$. Thus,
\begin{equation*}
\begin{array}{l}
\cr E[\emph{\textbf{x}}^{\emph{\textbf{X}}(t)}\emph{\textbf{y}}
^{\emph{\textbf{Y}}(t)}\emph{\textbf{z}}^{\emph{\textbf{Z}}(t)}\mid (\emph{\textbf{X}}(0),\emph{\textbf{Y}}(0),\emph{\textbf{Z}}(0))
=(\emph{\textbf{i}},\emph{\textbf{k}},\tilde{\emph{\textbf{k}}})]\\
\cr =E[\emph{\textbf{x}}^{\sum\limits_{a=1}^2\sum\limits_{k=1}^{i_a}
\emph{\textbf{X}}_{a,k}(t)}\emph{\textbf{y}}
^{\emph{\textbf{k}}+\sum\limits_{a=1}^2\sum\limits_{k=1}^{i_a}
\emph{\textbf{Y}}_{a,k}(t)}\emph{\textbf{z}}^{\tilde{\emph{\textbf{k}}}
+\sum\limits_{a=1}^2\sum\limits_{k=1}^{i_a}
\emph{\textbf{Z}}_{a,k}(t)}]\\
\cr=\emph{\textbf{y}}^{\emph{\textbf{k}}}\emph{\textbf{z}}
^{\tilde{\emph{\textbf{k}}}}E[\prod\limits_{k=1}^{i_1}\emph{\textbf{x}}
^{\emph{\textbf{X}}_{1,k}(t)}\prod\limits_{k=1}^{i_1}\emph{\textbf{y}}
^{\emph{\textbf{Y}}_{1,k}(t)}\prod\limits_{k=1}^{i_1}\emph{\textbf{z}}^{
\emph{\textbf{Z}}_{1,k}(t)}\cdot \prod\limits_{k=1}^{i_2}\emph{\textbf{x}}
^{\emph{\textbf{X}}_{2,k}(t)}\prod\limits_{k=1}^{i_2}\emph{\textbf{y}}
^{\emph{\textbf{Y}}_{2,k}(t)}\prod\limits_{k=1}^{i_2}\emph{\textbf{z}}^{
\emph{\textbf{Z}}_{2,k}(t)}]\\
\cr=\emph{\textbf{y}}^{\emph{\textbf{k}}}\emph{\textbf{z}}
^{\tilde{\emph{\textbf{k}}}}(E[\emph{\textbf{x}}
^{\emph{\textbf{X}}_{1,1}(t)}\emph{\textbf{y}}
^{\emph{\textbf{Y}}_{1,1}(t)}\emph{\textbf{z}}^{
\emph{\textbf{Z}}_{1,1}(t)}])^{i_1}\cdot (E[\emph{\textbf{x}}
^{\emph{\textbf{X}}_{2,1}(t)}\emph{\textbf{y}}
^{\emph{\textbf{Y}}_{2,1}(t)}\emph{\textbf{z}}^{
\emph{\textbf{Z}}_{2,1}(t)}])^{i_2}\\
\cr=\emph{\textbf{y}}^{\emph{\textbf{k}}}\emph{\textbf{z}}
^{\tilde{\emph{\textbf{k}}}}[\emph{\textbf{F}}(t,\emph{\textbf{x}},
\emph{\textbf{y}},\emph{\textbf{z}})]^{\emph{\textbf{i}}}.
\end{array}
\end{equation*}
The proof is complete. \hfill $\Box$
\end{proof}
\par
The functions $B_1(\emph{\textbf{x}},\emph{\textbf{y}})
+\bar{B}_1(\emph{\textbf{x}})$ and $B_2(\emph{\textbf{x}},\emph{\textbf{z}})
+\bar{B}_2(\emph{\textbf{x}})$ will play a significant role in the later discussion. The following theorem reveals their properties.
\par
\begin{Theorem} \label{th2.1}
{\rm(1)}\ For any $\textbf{y}\in [0,1)^{r_1},\textbf{z}\in [0,1)^{r_2}$,
\begin{eqnarray}\label{eq2.7}
\begin{cases}
B_1(\textbf{x},\textbf{y})+\bar{B}_1(\textbf{x})=0,\\
B_2(\textbf{x},\textbf{z})+\bar{B}_2(\textbf{x})=0
\end{cases}
\end{eqnarray}
possesses exact one root in ${[0,1]^2}$, denoted by $\textbf{q}(\textbf{y},\textbf{z}):=(q_1(\textbf{y},\textbf{z}),
q_2(\textbf{y},\textbf{z}))$. Moreover, $\textbf{q}(\textbf{y},\textbf{z})\leq \textbf{q}$, where $\textbf{q}=(q_1,q_2)$ is the minimal nonnegative solution of {\rm (\ref{eq2.a})} given in Lemma~\ref{le2.1}.
\par
{\rm(2)}\ $q_k(\textbf{y},\textbf{z})\in {C^\infty }([0,1)^{r_1+r_2})\ (k=1,2)$, and $q_k(\textbf{y},\textbf{z})$ can be expanded as a multivariate nonnegative Taylor series
\begin{eqnarray}\label{eq2.8}
q_k(\textbf{y},\textbf{z})=\sum\limits_{(\textbf{k,l})\in \mathbf{Z}_+^{r_1+r_2}}\beta^{(a)}_{\textbf{k,l}}\textbf{y}^{\textbf{k}}
\textbf{z}^{\textbf{l}},\quad (\textbf{y},\textbf{z})\in [0,1)^{r_1+r_2},\ \ k=1,2.
\end{eqnarray}
\end{Theorem}
\par
\begin{proof} Note that $B_1(\emph{\textbf{1}},\emph{\textbf{y}})+\bar{B}_1(\emph{\textbf{1}})<0$ and $B_2(\emph{\textbf{1}},\emph{\textbf{z}})+\bar{B}_2(\emph{\textbf{1}})<0$, by a similar argument as Lemma 2.8 in Li \& Wang~\cite{Li-W12}, we can prove that (\ref{eq2.7}) possesses exact one root in ${[0,1]^2}$. Note that
\begin{eqnarray*}
\begin{cases}
B_1(\emph{\textbf{x}},\emph{\textbf{y}})
+\bar{B}_1(\emph{\textbf{x}})\leq B_1(\emph{\textbf{x}}),\\
B_2(\emph{\textbf{x}},\emph{\textbf{z}})
+\bar{B}_2(\emph{\textbf{x}})\leq B_2(\emph{\textbf{x}}),
\end{cases}
\end{eqnarray*}
we further know that $\emph{\textbf{q}}(\emph{\textbf{y}},\emph{\textbf{z}})\leq \emph{\textbf{q}}$.
 \par
Next to prove (2). Integrating (\ref{eq2.5}) yields that for $k=1,2$,
\begin{eqnarray*}
&&\sum\limits_{(\emph{\textbf{j}},\emph{\textbf{k}},\tilde{\emph{\textbf{k}}})\in \mathbf{Z}_+^{2+r_1+r_2}}\tilde{p}_{_{(\emph{\textbf{e}}_k,\emph{\textbf{0}},
\tilde{\emph{\textbf{0}}}),
(\emph{\textbf{j}},\emph{\textbf{l}},\tilde{\emph{\textbf{l}}})}}(t)
\emph{\textbf{x}}^{\emph{\textbf{j}}}\emph{\textbf{y}}^{\emph{\textbf{l}}}
\emph{\textbf{z}}^{\tilde{\emph{\textbf{l}}}}-\emph{\textbf{x}}
^{\emph{\textbf{e}}_k}\nonumber\\
&=&[B_1(\emph{\textbf{x}},\emph{\textbf{y}})+\bar{B}
_1(\emph{\textbf{x}})]\int_0^t\frac{\partial F_{\emph{\textbf{e}}_k,\emph{\textbf{0}},\emph{\textbf{0}}}
(u,\emph{\textbf{x}},
\emph{\textbf{y}},\emph{\textbf{z}})}{{\partial {x_1}}}du+[B_2(\emph{\textbf{x}},\emph{\textbf{z}})+\bar{B}
_2(\emph{\textbf{x}})]\int_0^t\frac{\partial F_{\emph{\textbf{e}}_k,\emph{\textbf{0}},\emph{\textbf{0}}}(u,\emph{\textbf{x}},
\emph{\textbf{y}},\emph{\textbf{z}})}{{\partial {x_2}}}du.
\end{eqnarray*}
Since all the states $(\emph{\textbf{i}},\emph{\textbf{l}},\tilde{\emph{\textbf{l}}})$ with $|~\emph{\textbf{i}}~|>0$ are transient and all the states $(\emph{\textbf{0}},\emph{\textbf{l}},\tilde{\emph{\textbf{l}}})$ are absorbing, letting $\emph{\textbf{x}}=\emph{\textbf{q}}(\emph{\textbf{y}},\emph{\textbf{z}})$ in the above equality and then letting $t\rightarrow \infty$ yield that
\begin{eqnarray*}
q_k(\emph{\textbf{y}},
\emph{\textbf{z}})=\sum\limits_{(\emph{\textbf{k}},\tilde{\emph{\textbf{k}}})
\in \mathbf{Z}_+^{r_1+r_2}}\tilde{p}_{_{(\emph{\textbf{e}}_k,\emph{\textbf{0}},
\tilde{\emph{\textbf{0}}}),
(\emph{\textbf{0}},\emph{\textbf{l}},\tilde{\emph{\textbf{l}}})}}(+\infty)
\emph{\textbf{y}}^{\emph{\textbf{l}}}
\emph{\textbf{z}}^{\tilde{\emph{\textbf{l}}}}, \quad k=1,2.
\end{eqnarray*}
The proof is complete. \hfill $\Box$
\end{proof}

\par
\vspace{5mm}
\setcounter{section}{3}
\setcounter{equation}{0}
\setcounter{theorem}{0}
\setcounter{lemma}{0}
\setcounter{definition}{0}
\setcounter{corollary}{0}
\setcounter{remark}{0}
\setcounter{example}{0}
\noindent {\large \bf 3. Multiple birth property}
\vspace{3mm}
\par
Having prepared some preliminaries in the previous section, we now consider the multiple birth property of $2$-type Markov branching processes.
\par
We first give the following theorem which will play a key role in discussing the multiple birth property of $2$-type Markov branching processes.
\par
\begin{theorem}\label{th3.a}
Suppose that $\textbf{x}\in [0,1]^2, \textbf{y}\in [0,1)^{r_1}, [0,1)^{r_2}$.
\par
{\rm(1)}\ The differential equation
\begin{equation}\label{eq3.1}
\begin{cases}
\frac{\partial u_1}{\partial t}=B_1(\textbf{u},\textbf{y})
+\bar{B}_1(\textbf{u}),\\
\frac{\partial u_2}{\partial t}=B_2(\textbf{u},\textbf{z})
+\bar{B}_2(\textbf{u}),\\
\textbf{u}(0)=\textbf{x}
\end{cases}
\end{equation}
has unique solution $\textbf{u}(t)=\textbf{G}(t,\textbf{x},\textbf{y},\textbf{z})$, where \begin{equation*}
\textbf{u}(t)=(u_1(t),u_2(t)), \quad  \textbf{G}(t,\textbf{x},\textbf{y},\textbf{z})
=(g_1(t,\textbf{x},\textbf{y},\textbf{z}),
g_2(t,\textbf{x},\textbf{y},\textbf{z})).
\end{equation*}
\par
{\rm(2)}\ $\lim\limits_{t\rightarrow \infty}\textbf{G}(t,\textbf{x},\textbf{y},\textbf{z})=\textbf{q}(\textbf{y},
\textbf{z})$, where $\textbf{q}(\textbf{y},\textbf{z})$ is given in Theorem~\ref{th2.1}.
\par
\end{theorem}
\par
\begin{proof}
 We first prove (1). For fixed $(\emph{\textbf{y}},\emph{\textbf{z}})\in [0,1)^{r_1+r_2}$, denote
\begin{equation*}
\begin{cases}
H_1(\emph{\textbf{u}})=B_1(\emph{\textbf{u}},\emph{\textbf{y}}) +\bar{B}_1(\emph{\textbf{u}})-b_{\emph{\textbf{e}}_1}^{(1)}u_1,\\
H_2(\emph{\textbf{u}})=B_2(\emph{\textbf{u}},\emph{\textbf{z}})+ \bar{B}_2(\emph{\textbf{u}})-b_{\emph{\textbf{e}}_2}^{(2)}u_2.
\end{cases}
\end{equation*}
By the assumption (A-2), we know that $H_k(\emph{\textbf{u}})$ satisfies Lipchitz condition, i.e, there exists a constant $L$ such that for any $\emph{\textbf{u}}=(u_1,u_2), \tilde{\emph{\textbf{u}}}=(\tilde{u}_1,\tilde{u}_2)\in [0,1]^2$,
\begin{eqnarray*}
|H_k(\emph{\textbf{u}})-H_k(\tilde{\emph{\textbf{u}}})|\leq L\|\emph{\textbf{u}}-\tilde{\emph{\textbf{u}}}\|_1,\quad k=1,2,
\end{eqnarray*}
\par
For $\emph{\textbf{x}}\in [0,1]^2$, define $u_k^{(0)}(t)=x_ke^{b_{\emph{\textbf{e}}_k}^{(k)}t}\ (k=1,2)$ and
\begin{equation*} u_k^{(n)}(t)=e^{b_{\emph{\textbf{e}}_k}^{(k)}t}[x_k+\int_0^t
e^{-b_{\emph{\textbf{e}}_k}^{(k)}s}H_k(\emph{\textbf{u}}^{(n-1)}(s))ds],\quad n\geq 1,\ \ k=1,2.
\end{equation*}
\par
We can prove that
\begin{equation}\label{eq3.2}
0\leq u^{(n)}_k(t)\leq 1,\quad t\geq 0, n\geq 1,\ k=1,2
\end{equation}
and
\begin{eqnarray}\label{eq3.3}
\|\emph{\textbf{u}}^{(n+1)}(t)-\emph{\textbf{u}}^{(n)}(t)\|_1\leq \frac{M(2L)^n}{(n+1)!}t^{n+1},\quad t\geq 0,\ n\geq 1.
\end{eqnarray}
where $M:=|b^{(1)}_{\emph{\textbf{e}}_1}|+|b^{(2)}_{\emph{\textbf{e}}_2}|$. Indeed, it is obvious that $0\leq u^{(0)}_k(t)=x_ke^{b_{\emph{\textbf{e}}_k}^{(k)}t}\leq 1\ (k=1,2)$. Assume that
\begin{equation*}
0\leq u^{(n)}_k(t)\leq 1,\ \quad t\geq 0,\ k=1,2.
\end{equation*}
Then it is obvious that $u^{(n+1)}_k(t)\geq 0$ since $H_k(\emph{\textbf{u}})\geq 0$ for all $\emph{\textbf{u}}\in [0,1)^2$. On the other hand, for $k=1,2$,
\begin{eqnarray*}
u_k^{(n+1)}(t)&=&
e^{b_{\emph{\textbf{e}}_k}^{(k)}t}[x_k+\int_0^t e^{-b_{\emph{\textbf{e}}_k}^{(k)}s}H_k(\emph{\textbf{u}}^{(n)}(s))ds]\\
&\leq &e^{b_{\emph{\textbf{e}}_k}^{(k)}t}[x_k+\int_0^t e^{-b_{\emph{\textbf{e}}_k}^{(k)}s}H_k(\emph{\textbf{1}})ds]\\
&=&e^{b_{\emph{\textbf{e}}_k}^{(k)}t}[x_k-b_{\emph{\textbf{e}}_k}^{(k)}
\int_0^t e^{-b_{\emph{\textbf{e}}_k}^{(k)}s}ds]\\
&=&e^{b_{\emph{\textbf{e}}_k}^{(k)}t}[x_k+
e^{-b_{\emph{\textbf{e}}_k}^{(k)}t}-1]\\
&\leq &1.
\end{eqnarray*}
(\ref{eq3.2}) is proved. As for (\ref{eq3.3}), by the definition of $\emph{\textbf{u}}^{(n)}(t)$,
\begin{eqnarray*}
|u^{(n+1)}_k(t)-u^{(n)}_k(t)|&\leq& e^{b_{\emph{\textbf{e}}_k}^{(k)}t}\int_0^t
e^{-b_{\emph{\textbf{e}}_k}^{(k)}s}|H_k(\emph{\textbf{u}}^{(n)}(s))
-H_k(\emph{\textbf{u}}^{(n-1)}(s))|\ ds\\
&\leq& L\int_0^t\|\emph{\textbf{u}}^{(n)}(s)-\tilde{\emph{\textbf{u}}}^{(n-1)}(s)\|_1\ ds,\quad n\geq 1,\ k=1,2.
\end{eqnarray*}
Hence,
\begin{eqnarray}\label{eq3.4}
\|\emph{\textbf{u}}^{(n+1)}(t)-\emph{\textbf{u}}^{(n)}(t)\|_1&\leq& 2L\int_0^t\|\emph{\textbf{u}}^{(n)}(s)-\tilde{\emph{\textbf{u}}}^{(n-1)}(s)\|_1\ ds,\quad n\geq 1.
\end{eqnarray}
\par
Note that
\begin{eqnarray*}
|u^{(1)}_k(t)-u^{(0)}_k(t)|=e^{b_{\emph{\textbf{e}}_k}^{(k)}t}\int_0^t
e^{-b_{\emph{\textbf{e}}_k}^{(k)}s}H_k(\emph{\textbf{u}}^{(0)}(s))ds
\leq |b^{(k)}_{\emph{\textbf{e}}_k}|t,\quad k=1,2,
\end{eqnarray*}
we know that
\begin{eqnarray}\label{eq3.5}
\|\emph{\textbf{u}}_1(t)-\emph{\textbf{u}}_0(t)\|_1\leq Mt,
\end{eqnarray}
It follows from (\ref{eq3.4}), (\ref{eq3.5}) and mathematical induction that (\ref{eq3.3}) holds.
\par
Since
$$
u^{(n)}_k(t)=u^{(0)}_k(t)+\sum\limits_{j=1}^n(u^{(j)}_k(t)-u^{(j-1)}_k(t)),
\quad k=1,2,
$$
by (\ref{eq3.3}), we know that $u^{(n)}_k(t)\ (k=1,2)$
converges uniformly in any finite interval $[0,T]$. Therefore, $u_k(t):=\lim\limits_{n\rightarrow \infty}u^{(n)}_k(t)$ exists and it can be easily checked that $\emph{\textbf{u}}(t)=(u_1(t),u_2(t))$ is a solution of (\ref{eq3.1}). On the other hand, since $B_1(\emph{\textbf{u}},\emph{\textbf{y}}), \bar{B}_1(\emph{\textbf{u}})$, $B_2(\emph{\textbf{u}},\emph{\textbf{z}})$ and $\bar{B}_2(\emph{\textbf{u}})$ satisfy Lipchitz condition, by the differential equations theory, we know that (\ref{eq3.1}) has unique solution. The unique solution of (\ref{eq3.1}) is denoted by $\emph{\textbf{G}}(t,\emph{\textbf{x}},\emph{\textbf{y}},\emph{\textbf{z}})$.
\par
We now prove (2). For fixed $(\emph{\textbf{x}},\emph{\textbf{y}},\emph{\textbf{z}})\in [0,1]^2\times [0,1)^{r_1+r_2}$, denote
\begin{eqnarray*}
&&f_1(\emph{\textbf{u}}):=B_1(\emph{\textbf{u}},\emph{\textbf{y}})+\bar{B}_1
(\emph{\textbf{u}}),\\
&&f_2(\emph{\textbf{u}}):=B_2(\emph{\textbf{u}},\emph{\textbf{z}})+\bar{B}_2
(\emph{\textbf{u}}),\\
&& \emph{\textbf{G}}(t)=(g_1(t),g_2(t)):=\emph{\textbf{G}}(t,\emph{\textbf{x}},
\emph{\textbf{y}},\emph{\textbf{z}})
\end{eqnarray*}
for a moment.
\par
(a)\ Suppose that $f_1(\emph{\textbf{x}})\geq 0, f_2(\emph{\textbf{x}})\geq 0$. We prove that
\begin{eqnarray*}
\omega:=\inf\limits_{t\geq 0}\{\min(f_1(\emph{\textbf{G}}(t)),f_2(\emph{\textbf{G}}(t)))\}\geq 0.
\end{eqnarray*}
Indeed, suppose that $\omega<0$. Then by the continuity of $f_1, f_2$ and $\emph{\textbf{G}}(t)$, there exist $\tilde{t}<+\infty$ and $\delta>0$ such that
\begin{eqnarray}\label{eq3.6}
\min(f_1(\emph{\textbf{G}}(\tilde{t})),f_2(\emph{\textbf{G}}
(\tilde{t})))=0,\quad \min(f_1(\emph{\textbf{G}}(\tilde{t})+s),f_2(\emph{\textbf{G}}
(\tilde{t}+s)))<0,\ \ \forall s\in (0,\delta).
\end{eqnarray}
We can assume $f_1(\emph{\textbf{G}}(\tilde{t}))=0$ without loss of generality. If $f_2(\emph{\textbf{G}}(\tilde{t}))>0$, then there exists $\tilde{\delta}\in (0,\delta)$ such that
\begin{eqnarray*}
f_1(\emph{\textbf{G}}(\tilde{t}+s))<0, \quad f_2(\emph{\textbf{G}}(\tilde{t}+s))>0, \quad s\in (0,\tilde{\delta}),
\end{eqnarray*}
which, by (\ref{eq3.1}), implies that
\begin{eqnarray*}
g_1(\emph{\textbf{G}}(\tilde{t}+s))<g_1(\emph{\textbf{G}}(\tilde{t})),\quad
g_2(\emph{\textbf{G}}(\tilde{t}+s))>g_2(\emph{\textbf{G}}(\tilde{t})),\quad s\in (0,\tilde{\delta}).
\end{eqnarray*}
Therefore,
\begin{eqnarray}\label{eq3.7} f_1(g_1(\emph{\textbf{G}}(\tilde{t}+s)),g_2(\emph{\textbf{G}}
(\tilde{t})))\leq f_1(\emph{\textbf{G}}(\tilde{t}+s))<0, \quad s\in (0,\tilde{\delta}).
\end{eqnarray}
However, it is well-known that $u=g_1(\emph{\textbf{G}}
(\tilde{t}))$ is the unique root of $f_1(u,g_2(\emph{\textbf{G}}
(\tilde{t})))=0$ in $[0,1]$ with $f_1(u,g_2(\emph{\textbf{G}}
(\tilde{t})))>0$ for $u\in [0,g_1(\emph{\textbf{G}}
(\tilde{t})))$, which contradicts with (\ref{eq3.7}). Therefore,
\begin{eqnarray*}
f_1(\emph{\textbf{G}}(\tilde{t}))=0,\quad f_2(\emph{\textbf{G}}(\tilde{t}))=0.
\end{eqnarray*}
By Theorem~\ref{th2.1}, $\emph{\textbf{G}}(\tilde{t})=\emph{\textbf{q}}(\emph{\textbf{y}},\emph{\textbf{z}})$.
Hence, by (1), we know that $\emph{\textbf{G}}(t)=\emph{\textbf{q}}(\emph{\textbf{y}},\emph{\textbf{z}})$ for $t\geq \tilde{t}$. Thus,
\begin{eqnarray*}
f_1(\emph{\textbf{G}}(\tilde{t}+s))=f_2(\emph{\textbf{G}}(\tilde{t}+s))=0,\quad s\geq 0,
\end{eqnarray*}
which contradicts with (\ref{eq3.6}). Therefore, we have $\omega\geq 0$. Hence, $\emph{\textbf{G}}(t)$ is increasing in $t\geq 0$. By (\ref{eq3.1}),
\begin{eqnarray}\label{eq3.8a}
g_k(t)=e^{b_{\emph{\textbf{e}}_k}^{(k)}t}[x_k+\int_0^t
e^{-b_{\emph{\textbf{e}}_k}^{(k)}s}H_k(\emph{\textbf{G}}(s))ds],\quad k=1,2.
\end{eqnarray}
Letting $t\rightarrow \infty$ in the above equality yields
\begin{eqnarray*}
\begin{cases}
B_1(\lim\limits_{t\rightarrow\infty}\emph{\textbf{G}}(t),\emph{\textbf{y}})
+\bar{B}_1(\lim\limits_{t\rightarrow\infty}\emph{\textbf{G}}(t))=0\\
B_2(\lim\limits_{t\rightarrow\infty}\emph{\textbf{G}}(t),\emph{\textbf{z}})
+\bar{B}_2(\lim\limits_{t\rightarrow\infty}\emph{\textbf{G}}(t))=0.
\end{cases}
\end{eqnarray*}
Therefore,
\begin{eqnarray*}
\lim\limits_{t\rightarrow \infty}\emph{\textbf{G}}(t)=\emph{\textbf{q}}(\emph{\textbf{y}},
\emph{\textbf{z}}).
\end{eqnarray*}
\par
(b)\ Suppose that $f_1(\emph{\textbf{x}})\leq 0, f_2(\emph{\textbf{x}})\leq 0$. We can prove that
\begin{eqnarray*}
\omega:=\sup\limits_{t\geq 0}\{\min(f_1(\emph{\textbf{G}}(t)),f_2(\emph{\textbf{G}}(t)))\}\leq 0.
\end{eqnarray*}
By a similar argument as in (a), it can be proved that $\emph{\textbf{G}}(t)$ is decreasing in $t\geq 0$ and  \begin{eqnarray*}
\lim\limits_{t\rightarrow \infty}\emph{\textbf{G}}(t)=\emph{\textbf{q}}(\emph{\textbf{y}},
\emph{\textbf{z}}).
\end{eqnarray*}
\par
(c)\ Suppose that $f_1(\emph{\textbf{x}})\geq 0, f_2(\emph{\textbf{x}})< 0$. Let
\begin{eqnarray*}
\sigma=\inf\{t\geq 0: f_1(\emph{\textbf{G}}(t))\leq 0\ {\rm{or}}\ f_2(\emph{\textbf{G}}(t))\geq 0\}.
\end{eqnarray*}
If $\sigma<+\infty$, then $g_1(\emph{\textbf{G}}(t))$ is increasing and $g_2(\emph{\textbf{G}}(t))$ is decreasing in $[0,\sigma)$. It can be easily checked that $\emph{\textbf{G}}(\sigma+t)$ is the solution of (\ref{eq3.1}) with initial condition $\emph{\textbf{G}}(\sigma)$. Furthermore, we have that $f_1(\emph{\textbf{G}}(\sigma))\geq 0,\ f_2(\emph{\textbf{G}}(\sigma))=0$ or that $f_1(\emph{\textbf{G}}(\sigma))=0,\ f_2(\emph{\textbf{G}}(\sigma))<0$.
In the case that $f_1(\emph{\textbf{G}}(\sigma))\geq 0,\ f_2(\emph{\textbf{G}}(\sigma))=0$, by (a), we know that $g_1(\emph{\textbf{G}}(t))$ and $g_2(\emph{\textbf{G}}(t))$ are both
increasing in $t\in [\sigma,+\infty)$ and
\begin{eqnarray*}
\lim\limits_{t\rightarrow \infty}\emph{\textbf{G}}(t)=\emph{\textbf{q}}(\emph{\textbf{y}},
\emph{\textbf{z}}).
\end{eqnarray*}
while in the case that $f_1(\emph{\textbf{G}}(\sigma))=0,\ f_2(\emph{\textbf{G}}(\sigma))<0$, by (b), we know that $g_1(\emph{\textbf{G}}(t))$ and $g_2(\emph{\textbf{G}}(t))$ are both
decreasing in $t\in [\sigma,+\infty)$ and
\begin{eqnarray*}
\lim\limits_{t\rightarrow \infty}\emph{\textbf{G}}(t)=\emph{\textbf{q}}(\emph{\textbf{y}},
\emph{\textbf{z}}).
\end{eqnarray*}
\par
If $\sigma=+\infty$, then $g_1(\emph{\textbf{G}}(t))$ is increasing and $g_2(\emph{\textbf{G}}(t))$ is decreasing in $t\geq 0$. By (\ref{eq3.8a}), we still have
\begin{eqnarray*}
\lim\limits_{t\rightarrow \infty}\emph{\textbf{G}}(t)=\emph{\textbf{q}}(\emph{\textbf{y}},
\emph{\textbf{z}}).
\end{eqnarray*}
\par
(d)\ Suppose that $f_1(\emph{\textbf{x}})<0, f_2(\emph{\textbf{x}})\geq 0$. Let
\begin{eqnarray*}
\sigma=\inf\{t\geq 0: f_1(\emph{\textbf{G}}(t))\geq 0\ {\rm{or}}\ f_2(\emph{\textbf{G}}(t))\leq 0\}.
\end{eqnarray*}
A similar argument as in (c) yields the conclusion.
 The proof is complete.\hfill $\Box$
\end{proof}
\par
The following theorem gives the joint probability generating function of $(\emph{\textbf{Y}}(t),\emph{\textbf{Z}}(t))$.~\\
\par
\begin{theorem}\label{th3.1}
Suppose that $\{\textbf{X}(t):t \geq 0\} $ is a $2$-type Markov branching process with $\textbf{X}(0) = \textbf{e}_k, (k=1\ {\rm{or}}\ 2)$. $\textbf{G}(t,\textbf{x},\textbf{y},\textbf{z}) = (g_1(t,\textbf{x},\textbf{y},\textbf{z}),g_2(t,\textbf{x},\textbf{y},
\textbf{z}))$ is the unique solution of $(\ref{eq3.1})$. Then, the joint probability generating function of $(\textbf{Y}(t),\textbf{Z}(t))$ is given by
\begin{equation}\label{eq3.8}
E[\textbf{y}^{\textbf{Y}(t)}\textbf{z}^{\textbf{Z}(t)}\mid \textbf{X}(0)=\textbf{e}_k]=g_k(t,\textbf{1},\textbf{y},\textbf{z}),
\quad (\textbf{y},\textbf{z})\in [0,1)^{r_1+r_2},\ \ k=1,2.
\end{equation}
\par
In particular, the joint probability generating function of $\textbf{Y}(t)$ and $\textbf{Z}(t))$ are given by
\begin{equation}\label{eq3.9}
E[\textbf{y}^{\textbf{Y}(t)}\mid \textbf{X}(0)=\textbf{e}_k]=g_k(t,\textbf{1},\textbf{y},\textbf{1}),
\quad \textbf{y}\in [0,1)^{r_1},\ \ k=1,2.
\end{equation}
and
\begin{equation}\label{eq3.10}
E[\textbf{z}^{\textbf{Z}(t)}\mid \textbf{X}(0)=\textbf{e}_k]=g_k(t,\textbf{1},\textbf{1},\textbf{z}),
\quad \textbf{z}\in [0,1)^{r_2},\ \ k=1,2,
\end{equation}
respectively.
\end{theorem}
\par
\begin{proof}
Let $\tilde{P}(t)=(\tilde{p}_{_{(\emph{\textbf{i}},\emph{\textbf{k}},
\tilde{\emph{\textbf{k}}}),
(\emph{\textbf{j}},\emph{\textbf{l}},\tilde{\emph{\textbf{l}}})}}(t):
(\emph{\textbf{i}},\emph{\textbf{k}},\tilde{\emph{\textbf{k}}}),
(\emph{\textbf{j}},\emph{\textbf{l}},\tilde{\emph{\textbf{l}}})\in \mathbf{Z}_+^{2+r_1+r_2})$ be the transition probability of $(\emph{\textbf{X}}(t),\emph{\textbf{Y}}(t),\emph{\textbf{Z}}(t))$. We need to prove that for any fixed $(\emph{\textbf{x}},\emph{\textbf{y}},\emph{\textbf{z}}) \in [0,1]^{2+r_1+r_2}$,
\begin{equation}\label{eq3.11}
g_k(t,\emph{\textbf{x}},\emph{\textbf{y}},\emph{\textbf{z}})= F_k(t,\emph{\textbf{x}},\emph{\textbf{y}},\emph{\textbf{z}}),\quad k=1,2,
\end{equation}
where $F_k(t,\emph{\textbf{x}},\emph{\textbf{y}},\emph{\textbf{z}})\ (k=1,2)$ are given in Lemma~\ref{le2.3}.
It is sufficient to prove that for any $(\emph{\textbf{y}},\emph{\textbf{z}})\in [0,1)^{r_1+r_2}$,
\begin{equation*}
u_k(t,\emph{\textbf{x}}):=F_k(t,\emph{\textbf{x}},\emph{\textbf{y}},
\emph{\textbf{z}}),\quad k=1,2.
\end{equation*}
is a solution of (\ref{eq3.1}). Indeed, suppose $k=1$ without loss of generality, by Kolmogorov backward equation, for any $t \ge 0$, we have,
\begin{equation*}
\tilde{p}'_{_{(\emph{\textbf{e}}_1,\emph{\textbf{0}},
\tilde{\emph{\textbf{0}}}),
(\emph{\textbf{j}},\emph{\textbf{l}},\tilde{\emph{\textbf{l}}})}}(t)
=\sum\limits_{(\emph{\textbf{i}},\emph{\textbf{k}},\tilde{\emph{\textbf{k}}})
\in \mathbf{Z}_+^{2+r_1+r_2}}q_{_{(\emph{\textbf{e}}_1,\emph{\textbf{0}},
\tilde{\emph{\textbf{0}}}),
(\emph{\textbf{i}},\emph{\textbf{k}},\tilde{\emph{\textbf{k}}})}}
\tilde{p}_{_{(\emph{\textbf{i}},\emph{\textbf{k}},\tilde{\emph{\textbf{k}}}),
(\emph{\textbf{j}},\emph{\textbf{l}},\tilde{\emph{\textbf{l}}})}}(t).
\end{equation*}
Multiply $\emph{\textbf{x}}^{\emph{\textbf{j}}}\emph{\textbf{y}}^{\emph{\textbf{l}}}
\emph{\textbf{z}}^{\tilde{\emph{\textbf{l}}}}$ on both sides of the above equality and take summation over $(\emph{\textbf{j}},\emph{\textbf{l}},\tilde{\emph{\textbf{l}}})\in \mathbf{Z}_+^{2+r_1+r_2}$, we get
\begin{eqnarray*}
&&\sum\limits_{(\emph{\textbf{j}},\emph{\textbf{l}},\tilde{\emph{\textbf{l}}})
\in\mathbf{Z}_+^{2+r_1+r_2}}\tilde{p}'_{_{(\emph{\textbf{e}}_1,\emph{\textbf{0}},
\tilde{\emph{\textbf{0}}}),
(\emph{\textbf{j}},\emph{\textbf{l}},\tilde{\emph{\textbf{l}}})}}(t)
\emph{\textbf{x}}^{\emph{\textbf{j}}}\emph{\textbf{y}}^{\emph{\textbf{l}}}
\emph{\textbf{z}}^{\tilde{\emph{\textbf{l}}}}
= \sum\limits_{\emph{\textbf{i}}\in R_1}b^{(1)}_{\emph{\textbf{i}}}F_{_{\emph{\textbf{i}},\varepsilon_{
\emph{\textbf{i}}},
\tilde{\emph{\textbf{0}}}}}(t,\emph{\textbf{x}},\emph{\textbf{y}},
\emph{\textbf{z}})+\sum\limits_{\emph{\textbf{i}}\in R_1^c}b^{(1)}_{\emph{\textbf{i}}}F_{_{\emph{\textbf{i}},
\emph{\textbf{0}},
\tilde{\emph{\textbf{0}}}}}(t,\emph{\textbf{x}},\emph{\textbf{y}},
\emph{\textbf{z}})
\end{eqnarray*}
By (\ref{eq2.6}),
\begin{eqnarray*}
\frac{\partial F_1(t,\emph{\textbf{x}},\emph{\textbf{y}},\emph{\textbf{z}})}{\partial t}
=B_1(\emph{\textbf{F}}(t,\emph{\textbf{x}},\emph{\textbf{y}},\emph{\textbf{z}}),
\emph{\textbf{y}})+\bar{B}_1(\emph{\textbf{F}}(t,\emph{\textbf{x}},\emph{\textbf{y}},
\emph{\textbf{z}})).
\end{eqnarray*}
By a similar argument, we have
\begin{eqnarray*}
\frac{\partial F_2(t,\emph{\textbf{x}},\emph{\textbf{y}},\emph{\textbf{z}})}{\partial t}
=B_2(\emph{\textbf{F}}(t,\emph{\textbf{x}},\emph{\textbf{y}},\emph{\textbf{z}}),
\emph{\textbf{y}})+\bar{B}_2(\emph{\textbf{F}}(t,\emph{\textbf{x}},\emph{\textbf{y}},
\emph{\textbf{z}})).
\end{eqnarray*}
Note that $F_k(0,\emph{\textbf{x}},\emph{\textbf{y}},\emph{\textbf{z}})=x_k\ (k=1,2)$, we know that $u_k(t,\emph{\textbf{x}})=F_k(t,\emph{\textbf{x}},\emph{\textbf{y}},
\emph{\textbf{z}})\ (k=1,2)$ is a solution of (\ref{eq3.1}).
\par
Therefore, (\ref{eq3.11}) and hence (\ref{eq3.8}) holds. Finally, (\ref{eq3.9}) and (\ref{eq3.10}) follows directly from (\ref{eq3.8}).
The proof is complete.\hfill $\Box$
\end{proof}
\par
The following proposition presents the probability generating function of $(\emph{\textbf{Y}}(t),\emph{\textbf{Z}}(t))$ when the process $\emph{\textbf{t}}$ starts at $\emph{\textbf{X}}(0)=\emph{\textbf{i}}$.
\par
\begin{proposition}\label{pro3.1}
Suppose that $\{\textbf{X}(t):t \geq 0\} $ is a $2$-type Markov branching process with $\textbf{X}(0)=\textbf{i}$. Then,
\begin{equation}\label{eq3.12}
E[\textbf{y}^{\textbf{Y}(t)}\textbf{z}^{\textbf{Z}(t)}\mid \textbf{X}(0)=\textbf{i}]=[\textbf{G}(t,\textbf{1},\textbf{y},\textbf{z})]
^{\textbf{i}},
\quad (\textbf{y},\textbf{z})\in [0,1)^{r_1+r_2}.
\end{equation}
In particular,
\begin{equation}\label{eq3.13}
E[\textbf{y}^{\textbf{Y}(t)}\mid \textbf{X}(0)=\textbf{i}]=[\textbf{G}(t,\textbf{1},\textbf{y},\textbf{1})]
^{\textbf{i}},
\quad \textbf{y}\in [0,1)^{r_1}.
\end{equation}
and
\begin{equation}\label{eq3.14}
E[\textbf{z}^{\textbf{Z}(t)}\mid \textbf{X}(0)=\textbf{i}]=[\textbf{G}(t,\textbf{1},\textbf{1},\textbf{z})]
^{\textbf{i}},
\quad \textbf{z}\in [0,1)^{r_2}.
\end{equation}
\end{proposition}
\par
\begin{proof}
Since $E[\emph{\textbf{y}}^{\emph{\textbf{Y}}(t)}\emph{\textbf{z}}^{\emph{
\textbf{Z}}(t)}\mid \emph{\textbf{X}}(0)=\emph{\textbf{i}}]=F_{\emph{\textbf{i}},
\emph{\textbf{0}},
\tilde{\emph{\textbf{0}}}}(t,\emph{\textbf{1}},\emph{\textbf{y}},
\emph{\textbf{z}})$, by (\ref{eq2.6}) and Theorem~\ref{th3.1}, we immediately obtain (\ref{eq3.12}). (\ref{eq3.13}) and (\ref{eq3.14}) follows directly from (\ref{eq3.12}). The proof is complete. \hfill $\Box$
\end{proof}
\par
As direct consequences of Theorem~\ref{th3.1}, the following corollaries give the probability generating functions of the pure death number of type-$k$ individuals and twins-birth number of type-$k$ individuals.
\par
\begin{corollary}\label{cor3.1}
Suppose that $\{\textbf{X}(t):t \ge 0\}$ is a $2$-type Markov branching process with $\textbf{X}(0)=\textbf{e}_k\ (k=1,2)$, $Y(t)$ and $Z(t)$ are the pure death numbers of type-$1$ and type-$2$ individuals, respectively. Then,
\begin{equation}\label{eq3.15}
E[y^{Y(t)}z^{Z(t)}\mid \textbf{X}(0)=\textbf{e}_k]= g_k(t,y,z),\quad y,z\in [0,1),\ k=1,2.
\end{equation}
\par
In particular,
\begin{equation}\label{eq3.16}
E[y^{Y(t)}\mid \textbf{X}(0)=\textbf{e}_k]= g_k(t,y,1),\quad y\in [0,1),\ k=1,2
\end{equation}
and
\begin{equation}\label{eq3.17}
E[z^{Z(t)}\mid \textbf{X}(0)=\textbf{e}_k]= g_k(t,1,z),\quad z\in [0,1),\ k=1,2,
\end{equation}
where $(g_1(t,y,z),g_2(t,y,z))$ is the unique solution of the equation
\begin{eqnarray*}
\begin{cases}
\frac{\partial u_1}{\partial t}=B_1(u_1,u_2)- b_{\textbf{0}}^{(1)}(1-y)\\
\frac{\partial u_2}{\partial t}=B_2(u_1,u_2)- b_{\textbf{0}}^{(2)}(1-z)\\
u_1(0)=u_2(0)=1.
\end{cases}
\end{eqnarray*}
\end{corollary}
\par
\begin{proof}
Take $R_1=R_2=\{\emph{\textbf{0}}\}\subset \mathbf{Z}_+^2$. Then we have
\begin{eqnarray*}
&&B_1(\emph{\textbf{u}},y)+\bar{B}_1(\emph{\textbf{u}})=B_1(\emph{\textbf{u}})
-b^{(1)}_{\emph{\textbf{0}}}(1-y),\\
&&B_2(\emph{\textbf{u}},z)+\bar{B}_2(\emph{\textbf{u}})=B_2(\emph{\textbf{u}})
-b^{(1)}_{\emph{\textbf{0}}}(1-z).
\end{eqnarray*}
By Theorem~\ref{th3.1}, we immediately obtain (\ref{eq3.15}). (\ref{eq3.16}) and (\ref{eq3.17}) follows directly from (\ref{eq3.15}).
The proof is complete.\hfill $\Box$
\end{proof}
\par
\begin{corollary}\label{cor3.2}
Suppose that $\{\textbf{X}(t):t \ge 0\}$ is a $2$-type Markov branching process with $\textbf{X}(0)=\textbf{e}_k\ (k=1,2)$, $Y(t)$ is the $2\textbf{e}_1$-birth numbers of type-$1$ individuals and $Z(t)$ is the $2\textbf{e}_2$-birth numbers of type-$2$ individuals. Then,
\begin{equation*}
E[y^{Y(t)}z^{Z(t)}\mid \textbf{X}(0)=\textbf{e}_k]= g_k(t,y,z),\quad y,z\in [0,1),\ k=1,2.
\end{equation*}
\par
In particular,
\begin{equation*}
E[y^{Y(t)}\mid \textbf{X}(0)=\textbf{e}_k]= g_k(t,y,1),\quad y\in [0,1),\ k=1,2
\end{equation*}
and
\begin{equation*}
E[z^{Z(t)}\mid \textbf{X}(0)=\textbf{e}_k]= g_k(t,1,z),\quad z\in [0,1),\ k=1,2,
\end{equation*}
where $(g_1(t,y,z),g_2(t,y,z))$ is the unique solution of the equation
\begin{eqnarray*}
\begin{cases}
\frac{\partial u_1}{\partial t}=B_1(u_1,u_2)- b_{2\textbf{e}_1}^{(1)}(1-y)u_1^2\\
\frac{\partial u_2}{\partial t}=B_2(u_1,u_2)- b_{2\textbf{e}_2}^{(2)}(1-z)u_2^2\\
u_1(0)=u_2(0)=1.
\end{cases}
\end{eqnarray*}
\end{corollary}
\par
\begin{proof}
Take $R_1=\{2\emph{\textbf{e}}_1\}\subset \mathbf{Z}_+^2$ and $R_2=\{2\emph{\textbf{e}}_2\}\subset \mathbf{Z}_+^2$. Then we have
\begin{eqnarray*}
&&B_1(\emph{\textbf{u}},y)+\bar{B}_1(\emph{\textbf{u}})=B_1(\emph{\textbf{u}})
-b^{(1)}_{2\emph{\textbf{e}}_1}(1-y)u_1^2,\\
&&B_2(\emph{\textbf{u}},z)+\bar{B}_2(\emph{\textbf{u}})=B_2(\emph{\textbf{u}})
-b^{(2)}_{2\emph{\textbf{e}}_2}(1-z)u_2^2.
\end{eqnarray*}
By Theorem~\ref{th3.1}, we immediately obtain all the conclusions.
The proof is complete.\hfill $\Box$
\end{proof}
\par
Since $\emph{\textbf{0}}$ is the absorbing state of $\{\emph{\textbf{X}}(t):t \geq 0\}$, now we consider the multiple birth property until the extinction of the system. Let
\begin{equation*}
\tau=\inf\{t\geq 0: \emph{\textbf{X}}(t)=\emph{\textbf{0}}\}
\end{equation*}
be the extinction time of $\{\emph{\textbf{X}}(t):t \geq 0\}$.
\par
The following theorem gives the joint probability generating function of multi-birth number of individuals until the extinction of the system.
\par
\begin{theorem}\label{th3.3}
Suppose that $\{\textbf{X}(t):t \geq 0\}$ is a $2$-type Markov branching process with $\textbf{X}(0)=\textbf{e}_k\ (k=1,2)$.
\par
{\rm(1)} If $\rho(\textbf{1})\leq 0$, then the probability generating function of $(\textbf{Y}(\tau),\textbf{Z}(\tau))$ is given by
\begin{equation}\label{eq3.19}
E[\textbf{y}^{\textbf{Y}(\tau)}\textbf{z}^{\textbf{Z}(\tau)}\mid \textbf{X}(0)=\textbf{e}_k]=q_k(\textbf{y},\textbf{z}),\quad (\textbf{y},\textbf{z})\in [0,1)^{r_1+r_2},\ k=1,2,
\end{equation}
where $(q_1(\textbf{y},\textbf{z}),q_2(\textbf{y},\textbf{z}))$ is the unique solution of
\begin{eqnarray*}
\begin{cases}
B_1(\textbf{u},\textbf{y})+\bar{B}_1(\textbf{u})=0,\\
B_2(\textbf{u},\textbf{z})+\bar{B}_2(\textbf{u})=0.
\end{cases}
\end{eqnarray*}
\par
{\rm(2)}\ If $\rho(\textbf{1})>0$, then the probability generating function of $(\textbf{Y}(\tau),\textbf{Z}(\tau))$ conditioned on $\tau<\infty$ is given by
\begin{equation}
E[\textbf{y}^{\textbf{Y}(\tau)}\textbf{z}^{\textbf{Z}(\tau)}\mid \tau<\infty,\textbf{X}(0)=\textbf{e}_k]=\frac{q_k(\textbf{y},\textbf{z})}{q_k},\quad (\textbf{y},\textbf{z})\in [0,1)_+^{r_1+r_2},\ k=1,2.
\end{equation}
where $(q_1,q_2)$ is the minimal nonnegative solution of
\begin{eqnarray*}
\begin{cases}
B_1(\textbf{u})=0,\\
B_2(\textbf{u})=0.
\end{cases}
\end{eqnarray*}
\end{theorem}
\par
\begin{proof}
We first prove (1). It follows from Lemma~\ref{le2.3}(i) that for $k=1,2$ and any $(\emph{\textbf{x}},\emph{\textbf{y}},\emph{\textbf{z}}) \in [0,1]^2\times [0,1)^{r_1+r_2}$,
\begin{eqnarray*}
&&\sum\limits_{(\emph{\textbf{j}},\emph{\textbf{l}},
\tilde{\emph{\textbf{l}}})\in \mathbf{Z}_+^{2+r_1+r_2}}\tilde{p}_{_{(\emph{\textbf{e}}_k,\emph{\textbf{0}},
\tilde{\emph{\textbf{0}}}),
(\emph{\textbf{j}},\emph{\textbf{l}},\tilde{\emph{\textbf{l}}})}}(t)
\emph{\textbf{x}}^{\emph{\textbf{j}}}\emph{\textbf{y}}^{\emph{\textbf{l}}}
\emph{\textbf{z}}^{\tilde{\emph{\textbf{l}}}}-x_k\nonumber\\
&=&[B_1(\emph{\textbf{x}},\emph{\textbf{y}})+\bar{B}
_1(\emph{\textbf{x}})]\int_0^t\frac{\partial F_{\emph{\textbf{e}}_k,\emph{\textbf{0}},\tilde{\emph{\textbf{0}}}}
(s,\emph{\textbf{x}},\emph{\textbf{y}},
\emph{\textbf{z}})}{{\partial {x_1}}}ds+[B_2(\emph{\textbf{x}},\emph{\textbf{z}})+\bar{B}
_2(\emph{\textbf{x}})]\int_0^t\frac{\partial F_{\emph{\textbf{e}}_k,\emph{\textbf{0}},\tilde{\emph{\textbf{0}}}}(s,\emph{
\textbf{x}},\emph{\textbf{y}},
\emph{\textbf{z}})}{{\partial x_2}}ds.
\end{eqnarray*}
Letting $\emph{\textbf{x}}=\emph{\textbf{q}}(\emph{\textbf{y}},\emph{\textbf{z}})
=(q_1(\emph{\textbf{y}},\emph{\textbf{z}}),q_2(\emph{\textbf{y}},
\emph{\textbf{z}}))$
in the above equality and then letting $t\rightarrow \infty$ yield that
\begin{eqnarray*}
\sum\limits_{(\emph{\textbf{l}},
\tilde{\emph{\textbf{l}}})\in \mathbf{Z}_+^{r_1+r_2}}\tilde{p}_{_{(\emph{\textbf{e}}_k,\emph{\textbf{0}},
\tilde{\emph{\textbf{0}}}),
(\emph{\textbf{0}},\emph{\textbf{l}},\tilde{\emph{\textbf{l}}})}}(\infty)
\emph{\textbf{y}}^{\emph{\textbf{l}}}
\emph{\textbf{z}}^{\tilde{\emph{\textbf{l}}}}-q_k(\emph{\textbf{y}},
\emph{\textbf{z}})=0.
\end{eqnarray*}
\par
If $\rho(\emph{\textbf{1}})\leq 0$, then $q_k=P(\tau<\infty\mid \emph{\textbf{X}}(0)=\emph{\textbf{e}}_k)=1$. Therefore, noting that $(\emph{\textbf{0}},\emph{\textbf{l}},\tilde{\emph{\textbf{l}}})$ is absorbing state, we have
\begin{eqnarray*}
&&E[\emph{\textbf{y}}^{\emph{\textbf{Y}}(\tau)}\emph{\textbf{z}}
^{\emph{\textbf{Z}}
(\tau)}\mid \emph{\textbf{X}}(0)=\emph{\textbf{e}}_k]\\
&=&\sum\limits_{(\emph{\textbf{l}},
\tilde{\emph{\textbf{l}}})\in \mathbf{Z}_+^{r_1+r_2}}P((\emph{\textbf{Y}}(\tau),\emph{\textbf{Z}}(\tau))
=(\emph{\textbf{l}},\tilde{\emph{\textbf{l}}})\mid \emph{\textbf{X}}(0)=\emph{\textbf{e}}_k)
\emph{\textbf{y}}^{\emph{\textbf{l}}}
\emph{\textbf{z}}^{\tilde{\emph{\textbf{l}}}}\\
&=&\sum\limits_{(\emph{\textbf{l}},
\tilde{\emph{\textbf{l}}})\in \mathbf{Z}_+^{r_1+r_2}}\lim\limits_{t\rightarrow \infty}P((\emph{\textbf{Y}}(\tau),\emph{\textbf{Z}}(\tau))
=(\emph{\textbf{l}},\tilde{\emph{\textbf{l}}}), \tau<t\mid \emph{\textbf{X}}(0)=\emph{\textbf{e}}_k)
\emph{\textbf{y}}^{\emph{\textbf{l}}}
\emph{\textbf{z}}^{\tilde{\emph{\textbf{l}}}}\\
&=&\sum\limits_{(\emph{\textbf{l}},
\tilde{\emph{\textbf{l}}})\in \mathbf{Z}_+^{r_1+r_2}}\lim\limits_{t\rightarrow \infty}P((\emph{\textbf{Y}}(t),\emph{\textbf{Z}}(t))
=(\emph{\textbf{l}},\tilde{\emph{\textbf{l}}}), \tau<t\mid \emph{\textbf{X}}(0)=\emph{\textbf{e}}_k)
\emph{\textbf{y}}^{\emph{\textbf{l}}}
\emph{\textbf{z}}^{\tilde{\emph{\textbf{l}}}}\\
&=&\sum\limits_{(\emph{\textbf{l}},
\tilde{\emph{\textbf{l}}})\in \mathbf{Z}_+^{r_1+r_2}}\lim\limits_{t\rightarrow \infty}\tilde{p}_{_{(\emph{\textbf{e}}_k,\emph{\textbf{0}},
\tilde{\emph{\textbf{0}}}),
(\emph{\textbf{0}},\emph{\textbf{l}},\tilde{\emph{\textbf{l}}})}}(t)
\emph{\textbf{y}}^{\emph{\textbf{l}}}
\emph{\textbf{z}}^{\tilde{\emph{\textbf{l}}}}\\
&=&\sum\limits_{(\emph{\textbf{l}},
\tilde{\emph{\textbf{l}}})\in \mathbf{Z}_+^{r_1+r_2}}\tilde{p}_{_{(\emph{\textbf{e}}_k,\emph{\textbf{0}},
\tilde{\emph{\textbf{0}}}),
(\emph{\textbf{0}},\emph{\textbf{l}},\tilde{\emph{\textbf{l}}})}}(\infty)
\emph{\textbf{y}}^{\emph{\textbf{l}}}
\emph{\textbf{z}}^{\tilde{\emph{\textbf{l}}}}\\
&=& q_k(\emph{\textbf{y}},
\emph{\textbf{z}}).
\end{eqnarray*}
\par
(i) is proved.
\par
Next we prove (ii). If $\rho(\emph{\textbf{1}})\leq 0$, then $q_k=P(\tau<\infty\mid \emph{\textbf{X}}(0)=\emph{\textbf{e}}_k)<1$. Therefore, similarly as the above argument, we have
\begin{eqnarray*}
&&E[\emph{\textbf{y}}^{\emph{\textbf{Y}}(\tau)}\emph{\textbf{z}}
^{\emph{\textbf{Z}}
(\tau)}\mid \tau<\infty, \emph{\textbf{X}}(0)=\emph{\textbf{e}}_k]\\
&=&q_k^{-1}\sum\limits_{(\emph{\textbf{l}},
\tilde{\emph{\textbf{l}}})\in \mathbf{Z}_+^{r_1+r_2}}P((\emph{\textbf{Y}}(\tau),\emph{\textbf{Z}}(\tau))
=(\emph{\textbf{l}},\tilde{\emph{\textbf{l}}}),\tau<\infty\mid \emph{\textbf{X}}(0)=\emph{\textbf{e}}_k)
\emph{\textbf{y}}^{\emph{\textbf{l}}}
\emph{\textbf{z}}^{\tilde{\emph{\textbf{l}}}}\\
&=&q_k^{-1}\sum\limits_{(\emph{\textbf{l}},
\tilde{\emph{\textbf{l}}})\in \mathbf{Z}_+^{r_1+r_2}}\lim\limits_{t\rightarrow \infty}P((\emph{\textbf{Y}}(\tau),\emph{\textbf{Z}}(\tau))
=(\emph{\textbf{l}},\tilde{\emph{\textbf{l}}}), \tau<t\mid \emph{\textbf{X}}(0)=\emph{\textbf{e}}_k)
\emph{\textbf{y}}^{\emph{\textbf{l}}}
\emph{\textbf{z}}^{\tilde{\emph{\textbf{l}}}}\\
&=& \frac{q_k(\emph{\textbf{y}},\emph{\textbf{z}})}{q_k}.
\end{eqnarray*}
The proof is complete.\hfill $\Box$
\end{proof}
\par
By Theorem~\ref{th3.3}, we immediately obtain the following corollaries which gives the probability generating functions of the pure death number of type-$k$ individuals until the extinction of the system and twins-birth number of type-$k$ individuals until the extinction of the system.
\par
\begin{corollary}\label{cor3.3}
Suppose that $\{\textbf{X}(t):t \ge 0\}$ is a $2$-type Markov branching process with $\textbf{X}(0)=\textbf{e}_k\ (k=1,2)$, $Y(t)$ and $Z(t)$ are the pure death numbers of type-$1$ and type-$2$ individuals, respectively. If $\rho(\textbf{1})\leq 0$, then
\begin{equation*}
E[y^{Y(\tau)}z^{Z(\tau)}\mid \textbf{X}(0)=\textbf{e}_k]= q_k(y,z),\quad y,z\in [0,1),\ k=1,2.
\end{equation*}
\par
If $\rho(\textbf{1})>0$, then
\begin{equation*}
E[y^{Y(\tau)}z^{Z(\tau)}\mid \tau<\infty,\textbf{X}(0)=\textbf{e}_k]= \frac{q_k(y,z)}{q_k},\quad y,z\in [0,1),\ k=1,2.
\end{equation*}
where $(q_1(y,z),q_2(y,z))$ is the unique solution of the equation
\begin{eqnarray*}
\begin{cases}
B_1(u_1,u_2)- b_{\textbf{0}}^{(1)}(1-y)=0\\
B_2(u_1,u_2)- b_{\textbf{0}}^{(2)}(1-z)=0.
\end{cases}
\end{eqnarray*}
\end{corollary}
\par
\begin{proof}
Note $R_1=R_2=\{\emph{\textbf{0}}\}$, we immediately get the conclusions.
\hfill $\Box$
\end{proof}
\par
\begin{corollary}\label{cor3.4}
Suppose that $\{\textbf{X}(t):t \ge 0\}$ is a $2$-type Markov branching process with $\textbf{X}(0)=\textbf{e}_k\ (k=1,2)$, $Y(t)$ is the $2\textbf{e}_1$-birth numbers of type-$1$ individuals and $Z(t)$ is the $2\textbf{e}_2$-birth numbers of type-$2$ individuals. If $\rho(\textbf{1})\leq 0$, then
\begin{equation*}
E[y^{Y(\tau)}z^{Z(\tau)}\mid \textbf{X}(0)=\textbf{e}_k]= q_k(y,z),\quad y,z\in [0,1),\ k=1,2.
\end{equation*}
\par
If $\rho(\textbf{1})>0$, then
\begin{equation*}
E[y^{Y(\tau)}z^{Z(\tau)}\mid \tau<\infty,\textbf{X}(0)=\textbf{e}_k]= \frac{q_k(y,z)}{q_k},\quad y,z\in [0,1),\ k=1,2.
\end{equation*}
where $(q_1(y,z),q_2(y,z))$ is the unique solution of the equation
\begin{eqnarray*}
\begin{cases}
B_1(u_1,u_2)- b_{2\textbf{e}_1}^{(1)}(1-y)u_1^2=0\\
B_2(u_1,u_2)- b_{2\textbf{e}_2}^{(2)}(1-z)u_2^2=0.
\end{cases}
\end{eqnarray*}
\end{corollary}
\par
\begin{proof}
Note $R_1=\{2\emph{\textbf{e}}_1\}$ and $R_2=\{2\emph{\textbf{e}}_2\}$, we immediately get the conclusions.
\hfill $\Box$
\end{proof}
\par
Finally, we give an example to illustrate the main results obtained.
\par
\begin{example}
Suppose that $\{\emph{\textbf{X}}(t):t\geq 0\}$ is a $2$-type birth-death branching process with
$$
B_1(\emph{\textbf{x}})=p-x_1+qx_2^2,\quad B_2(\emph{\textbf{x}})=\alpha-x_2+\beta x_1,
$$
where $p,\ \alpha \in(0,1),\ q=1-p,\ \beta=1-\alpha$. $Y(t)$ is the pure death number of type-$1$ individuals until time $t$ and $Z(t)$ is the pure death number of type-$2$ individuals until time $t$. By Corollary~\ref{cor3.1}, we know that
\begin{eqnarray*}
E[y^{Y(t)}z^{Z(t)}\mid \emph{\textbf{X}}(0)=\emph{\textbf{e}}_k]=
\begin{cases}
u(t,y,z),\ & k=1,\\
v(t,y,z),\ & k=2,
\end{cases},\quad y,z\in [0,1),
\end{eqnarray*}
where $(u(t,y,z),v(t,y,z))$ is the unique solution of
\begin{eqnarray*}
\begin{cases}
\frac{\partial u}{\partial t}=qv^2-u+py\\
\frac{\partial v}{\partial t}=\beta u-v+\alpha z\\
u(0)=v(0)=1.
\end{cases}
\end{eqnarray*}
\par
It is easy to see that the maximum eigenvalue of $(B_{ij}(\emph{\textbf{1}}):i,j=1,2)$ is $\rho(\emph{\textbf{1}})=\sqrt{2q\beta}-1$. For $y,z\in [0,1)$, solving the equation
\begin{eqnarray*}
\begin{cases}
qv^2-u+py=0,\\
\beta u-v+\alpha z=0,
\end{cases}
\end{eqnarray*}
yields that
\begin{eqnarray*}
&&u=u(y,z)=\frac{1}{2q\beta^2}[1-\sqrt{1-4q\beta(p\beta y+\alpha z)}]-\frac{\alpha z}{\beta},\\
&&v=v(y,z)=\frac{1}{2q\beta}[1-\sqrt{1-4q\beta(p\beta y+\alpha z)}].
\end{eqnarray*}
By Corollary~\ref{cor3.3}, if $2q\beta\leq 1$, then
\begin{eqnarray*}
E[y^{Y(\tau)}z^{Z(\tau)}\mid\emph{\textbf{X}}(0)=\emph{\textbf{e}}_1]
=\frac{1-\sqrt{1-4q\beta(p\beta y+\alpha z)}-2q\beta \alpha z}{2q\beta^2},\quad y,z\in [0,1),
\end{eqnarray*}
\begin{eqnarray*}
E[y^{Y(\tau)}z^{Z(\tau)}\mid\emph{\textbf{X}}(0)=\emph{\textbf{e}}_2]
=\frac{1-\sqrt{1-4q\beta(p\beta y+\alpha z)}}{2q\beta},\quad y,z\in [0,1),
\end{eqnarray*}
\par
If $2q\beta>1$, then
\begin{eqnarray*}
E[y^{Y(\tau)}z^{Z(\tau)}\mid\emph{\textbf{X}}(0)=\emph{\textbf{e}}_1]
=\frac{1-\sqrt{1-4q\beta(p\beta y+\alpha z)}-2q\beta\alpha z}{2(1-2q\beta+q\beta^2)},\quad y,z\in [0,1),
\end{eqnarray*}
\begin{eqnarray*}
E[y^{Y(\tau)}z^{Z(\tau)}\mid\emph{\textbf{X}}(0)=\emph{\textbf{e}}_2]
=\frac{1-\sqrt{1-4q\beta(p\beta y+\alpha z)}}{2(1-q\beta)},\quad y,z\in [0,1).
\end{eqnarray*}
\end{example}

\section*{Acknowledgement}
\par
 This work is substantially supported by the National Natural Sciences Foundations of China (No. 11771452, No. 11971486).
\par
\noindent
\begin{center}
{\bf Declarations}
\end{center}
\hspace*{0.1cm}
\par
\noindent
{\bf Ethics Approval}\ Not applicable.
\par
\noindent
{\bf Funding}\ Funding provided by the National Natural Science Foundation of China (No. 11771452, No. 11971486).
\par
\noindent
{\bf Availability of data and materials}\ Not applicable.
\par
\noindent
{\bf Conflict of Interests}\ The authors declare that they have no conflict of interest.
\par

\end{document}